\documentclass[12pt]{article}

\usepackage{amstext    }
\usepackage{amsthm    }
\usepackage{a4}
\usepackage[mathscr]{eucal}
\usepackage{mathrsfs}

\usepackage{amsmath}
\usepackage{amssymb}
\usepackage{amscd}

\newtheorem{theorem}{Theorem}[section]

\newtheorem{proposition}[theorem]{Proposition}
\newtheorem{corollary}[theorem]{Corollary}
\newtheorem{lemma}[theorem]{Lemma}
\newtheorem{remark}[theorem]{Remark}

\newcommand{\cali}[1]{\mathscr{#1}}

\newcommand{\rank}{{\rm rank}}

\newcommand{\volume}{{\rm volume}}

\newcommand{\supp}{{\rm supp}}

\newcommand{\dbar}{{\overline\partial}}
\newcommand{\ddbar}{{\partial\overline\partial}}

\newcommand{\Bc}{\cali{B}}

\newcommand{\Gc}{\cali{G}}

\newcommand{\Kc}{\cali{K}}

\newcommand{\Qc}{\cali{Q}}

\newcommand{\C}{\mathbb{C}}
\newcommand{\Q}{\mathbb{Q}}

\newcommand{\N}{\mathbb{N}}
\newcommand{\Z}{\mathbb{Z}}
\newcommand{\R}{\mathbb{R}}

\renewcommand{\P}{\mathbb{P}}

\title{{\bf Dynamics of automorphisms on compact K\"ahler manifolds}}
\author{Henry De Th\'elin and Tien-Cuong Dinh}

\begin{document}
\maketitle

\begin{abstract}
We study holomorphic automorphisms on compact K\"ahler manifolds having simple actions on the Hodge cohomology ring. We show for such automorphisms
that the main dynamical Green currents admit  complex 
laminar structures (woven currents) and the Green measure is the unique invariant probability measure of
maximal entropy.
\end{abstract}

\noindent
{\bf Key-words:} holomorphic automorphism, Green current, laminar current, Green measure, entropy.

\medskip
\noindent
{\bf AMS Classification:} 32U40, 32H50.

\section{Introduction}

Let $(X,\omega)$ be a compact K\"ahler manifold of dimension $k$ where $\omega$ 
denotes a K\"ahler form on $X$. Let $f: X \to X$ be a holomorphic automorphism on $X$. 
The first part of this paper deals with some geometric property of dynamical Green currents  associated with $f$. These currents were constructed by Sibony and the second author in \cite{DS2}. 
We will show that they are woven. Roughly speaking, woven currents admit a complex laminar structure in the sense that they are averages of currents of integration on complex manifolds (see Section \ref{lamination} for precise definition).  The property is fundamental in the dynamical study of $f$ using a geometric method (see
Bedford-Lyubich-Smillie \cite{BLS}).

When $X$ is a projective manifold, the problem was solved in \cite{Di} using that $X$ admits many submanifolds of any dimension, see also Cantat \cite{C} for the case of dimension $k=2$.
For a general K\"ahler manifold, the approach breaks down and it is necessary to use another technique. 
Our approach here uses in particular a recent result by the first author which gives a criterion for a current to be woven \cite{Det6}. The criterion is valid in the local setting and allows us to work with general K\"ahler manifolds. 

The second part of the paper deals with the uniqueness of invariant measure of maximal entropy. This property  says that somehow the measure describes totally the most chaotic part of the dynamical system.
Our approach follows partially the method developed in \cite{BLS, DT2}. One of our 
contribution here concerns some new
results on the equidistribution  towards the Green currents of $f$. 
A priori, these currents are not intersections of currents of bidegree 
$(1,1)$ and this is a source of difficulties that we have to overcome in this work.

Recall that {\it  the dynamical degree} $d_p$ of order $p$ of 
$f$ is the spectral radius of the pull-back operator 
$f^*$ acting on the Hodge cohomology group $H^{p,p}(X,\C)$ for $0\leq p\leq k$. We have  $d_0=d_k=1$. 
An inequality due to 
Khovanskii, Teissier and Gromov \cite{Gromov} implies that the function $p \mapsto \log d_p$ is concave on $0 \leq p \leq k$. In particular, there are integers $s$ and $s'$ with $0 \leq  s \leq s' \leq k$ such that
$$1=d_0 < \cdots < d_s= \cdots = d_{s'} > \cdots > d_k=1.$$

In what follows, we assume that the action of $f^*$ on the Hodge cohomology ring
$\oplus H^*(X,\C)$ admits a unique eigenvalue of modulus $d_s$ which is moreover a simple eigenvalue. 
We say that the action of $f^*$ on Hodge cohomology group is {\it simple}. 
This property is equivalent to the fact that the sequence of linear operators  $d_s^{-n}(f^n)^*$ 
on $\oplus H^*(X,\C)$
converges to a rank 1 operator. The eigenvalue of maximal modulus is then equal to $d_s$ and we have $s=s'$, see \cite{DS2, DS1, G} for details.

An example of maps which do not satisfy the above property is the automorphism $(y,z)\mapsto (y,g(y))$ on the product $Y\times Z$ of two compact K\"ahler manifolds where $g$ is an automorphism on $Z$. So, our hypothesis on $f^*$ somehow insures  the lack of neutral direction in the dynamical system.
Under this natural condition, it was shown in
\cite{DS2} that the sequences $d_s^{-n}(f^n)^* \omega^s$ and $d_s^{-n}(f^n)_* \omega^{k-s}$ 
converge to positive closed currents $T^+$ and $T^-$ respectively. We call them {\it the main Green currents} associated with $f$. We will not consider in this work the Green currents of other bidegree.
Our first main result is the following.  

\begin{theorem}{\label{th1}}
Let $f$ be a holomorphic automorphism on a compact K\"ahler manifold $X$. 
Assume that the action of the pull-back operator $f^*$ on Hodge cohomology is simple. 
Then the main Green currents $T^\pm$ of $f$ are woven. If $T^+$ (resp. $T^-$) is of bidegree $(1,1)$, then it is laminar.
\end{theorem}

The intersection $\mu:=T^+\wedge T^-$ of $T^+$ and $T^-$ is well-defined and is a non-zero invariant positive measure. Multiplying $\omega$ with a constant allows us to assume that $\mu$ is a measure of probability. 
We call it {\it the Green measure} of $f$. The well-known variational principle says that the entropy of $\mu$ is bounded from above by the topological entropy of $f$ which is equal to $\log d_s$ according to results by 
Gromov and Yomdin \cite{Gromov2, Yomdin}. It was shown in \cite{DS1} that $\mu$ is a measure of maximal entropy, i.e. its entropy is equal to $\log d_s$. Here is our second main result.

\begin{theorem}{\label{th2}}
Let $f$ be as in Theorem \ref{th1}. Then the Green measure of $f$ is the unique invariant probability measure 
of maximal entropy.
\end{theorem}

\noindent
{\bf Acknowledgement.} This paper was partially written during the visit of the second author at the Max-Planck institute. He would like to thank this organization for its support and its hospitality.

\section{Laminar and woven currents}   \label{lamination}

In this section, we introduce the notion of laminar/woven current and give some criteria for currents to be laminar/woven. 

Let $X$ be a complex manifold of dimension $k$ and consider a Hermitian metric $\omega$ on $X$. 
Consider a connected complex manifold $Z$ of dimension $k-p$ with $1\leq p\leq k-1$ and a holomorphic map 
$g:Z\to X$. Assume that 
$$\int_Z g^*(\omega^{k-p}_{|K})<+\infty$$
for any compact subset $K$ of $X$. According to the Wirtinger's theorem, 
this condition means that the $2(k-p)$-dimensional volume of $g(Z)$ counted with multiplicity is locally finite in $X$. So, $g_*[Z]$ defines a positive $(p,p)$-current in $X$ which is not closed in general. Here, $[Z]$ is the current of integration on $Z$. 

We refer the reader 
to Demailly \cite{Demailly} and Federer \cite{Federer} for basic theory on currents. 
The above integral on $Z$ is the mass of $g_*[Z]$ on $K$. 
Recall that {\it the mass} on a Borel set $B\subset X$ of a positive $(p,p)$-current $S$ on $X$ is defined by
$$\|S \|_B:=\big\langle S,\omega^{k-p}_{| B}\big\rangle.$$

Denote by $\Bc_p(X)$ the set of all the currents of the form $g_*[Z]$ as above. Consider a positive measure $\nu$ on $\Bc_p(X)$ such that for any compact subset $K$ of $X$, we have 
$$\int \|S\|_K d\nu(S)<\infty.$$
Under this condition, we can define a positive $(p,p)$-current $T$ by
$$T:=\int S d\nu(S),$$
that is, 
$$\langle T,\alpha\rangle=\int \langle S,\alpha\rangle d\nu(S)$$
for all test $(k-p,k-p)$-form $\alpha$ on $X$.
The condition on $\nu$ insures that the last integral is meaningful. 
Such a current $T$ is called {\it woven}, see \cite{Di}.  If moreover, for $\nu$-almost all currents $S=g_*[Z]$ and $S'=g'_*[Z']$, the intersection $g(Z)\cap g'(Z')$ is empty or of maximal dimension $k-p$, we say that $T$ is {\it laminar}, see \cite{BLS}.

It is not difficult to see that $g_*[Z]$ can be written as a sum of currents of the same type with small support. Therefore, a current is woven if and only if it is locally woven. 
We will now give an equivalent definition of woven current which gives a simpler geometric structure  and is closer to the notion of laminar current introduced in \cite{BLS}. For simplicity, we assume that $X$ is an open set in $\C^k$. Otherwise, we can write $X$ as a finite or countable union $\overline X_i$ where $X_i$ are disjoint charts such that $T$ has no mass on $\partial X_i$ and  Proposition \ref{woven_def} below will give a description of $T$ in each  $X_i$. 

Consider a coordinates system $(z_1,\ldots,z_k)$ of $\C^k$ and the associated real coordinates systems
$(x_1,\ldots,x_{2k})$ of $\C^k \simeq \R^{2k}$ with $z_j=x_j+ix_{j+k}$. 
Divide $\C^k$ into cubes of size $r$ using the real hyperplanes $\{x_i=mr\}$ with $m\in\Z$ and $1\leq i\leq 2k$.
Such an $r$-cube is called {\it basic $r$-cube} and their union is denoted by $\Qc_r$.
We call this {\it a division} of $\C^k$ into $r$-cubes.
We can choose the coordinates system so that for all $r\in \Q$, the current $T$ has no mass on the complementary of $\Qc_r$, i.e. in the union of the hyperplanes $\{x_i=mr\}$.  

Consider an irreducible submanifold $\Gamma$ of codimension $p$ in a basic $r$-cube $D$. We say that $\Gamma$ is {\it nice} if it satisfies the following properties:
\begin{enumerate}
\item[(1)]  There is an irreducible submanifold $\Gamma'$ of the basic $3r$-cube $D'$ containing $\overline D$ such that $\Gamma$ is an open subset of $\Gamma'$;
\item[(2)] The manifold $\Gamma'$ is the graph of a holomorphic map over its tangent space at each point. 
\end{enumerate}
It is not difficult to see that the volume of such a manifold $\Gamma$ is bounded by a universal constant times $r^{2(k-p)}$. 

Denote by $\Bc_p(r)$ the set of all the nice manifolds $\Gamma$ of codimension $p$. Let $\nu$ be a positive measure on $\Bc_p(r)$ such that 
$$\int \volume(\Gamma) d\nu(\Gamma)<\infty.$$
We say that the associated current $S:=\int [\Gamma] d\nu(\Gamma)$ is {\it a nice woven current}, where $[\Gamma]$ denotes the current of integration on $\Gamma$. 

\begin{proposition} \label{woven_def}
Let $T$ be a woven $(p,p)$-current on an open subset $X$ of $\C^k$ as above. Then there are nice woven currents 
$T_i$ on $\Qc_{2^{-i}}$ which vanish outside $X$ and such that 
$$T=\sum_{i=0}^\infty T_i.$$
\end{proposition}
\proof
Consider first the case where $T$ is equal to the current $g_*[Z]$ as above. If $\dim g(Z)<k-p$, then $g_*[Z]=0$ as $(p,p)$-currents. So, we can assume that the dimension of $g(Z)$ is $k-p$. Denote by $\Sigma$ the set of points in $Z$ where the differential of $g$ is not of maximal rank. This is a proper analytic subset of $Z$. Therefore, its $2(k-p)$-dimensional volume vanishes. It follows that 
$[Z\setminus\Sigma]=[Z]$ and $g_*[Z\setminus \Sigma]=g_*[Z]$. 

Define $\Gc_i:=g^{-1}(\C^k\setminus \Qc_{2^{-i}})$. This is an increasing sequence of closed subsets  
in $Z$. Their volumes are equal to 0 because $g_*[Z]$ has no mass on the complementary of $\Qc_{2^{-i}}$. 
Denote by $W_0$ the union of connected components of $Z\setminus \Gc_0$ which are sent injectively by $g$ to nice manifolds in basic $1$-cubes. Define also by induction for $i\geq 1$,
$Z_i:=Z_{i-1}\setminus \overline W_{i-1}$ and 
$W_i$ the union of connected components of $Z_i\setminus \Gc_i$ which are sent injectively by $g$ to a nice manifolds in basic $2^{-i}$-cubes.

Observe that for any $a\in Z\setminus\Sigma$, if $r$ is small enough, there is a neighbourhood $U$ of $a$ such that $U$ is sent injectively to a nice submanifold of a basic $r$-cube. Therefore, the union of $W_i$ contains the complementary of $\cup_{i\geq 0} \Gc_i\cup\Sigma$. Since the last set has zero volume, we obtain that
$T=\sum T_i$ where $T_i:=g_*(W_i)$. The currents $T_i$ are nice woven in $\Qc_{2^{-i}}$. So, the proposition is true for $g_*[Z]$. 

In general, write $T=\int S d\nu(S)$ where $\nu$ is a positive measure on $\Bc_p(X)$. Since $T$ has no mass on the $\C^k\setminus \Qc_{2^{-i}}$, the same property holds for $\nu$-almost every $S$. Consider the decomposition $S=\sum S_i$ into nice currents obtained as above and define $T_i:=\int S_i d\nu(S)$.
It is clear that $T_i$ are nice woven currents and that $T=\sum T_i$.
\endproof

Woven currents appear naturally as certain limits of submanifolds. We give now a criterium for a sequence of manifolds to be convergent towards a woven current. It was obtained by the second author in the case of projective manifolds \cite{Di} and generalized by the first author to the local setting \cite{Det6}.   

Consider a sequence $M_n$ of (smooth) submanifolds of dimension $k-p$ in an open set $X$ of $\C^k$. 
Denote by $G(k-p,k)$ the set of complex subspaces of dimension $k-p$ in $\C^k$ through the origin $0\in\C^k$. 
This is a complex Grassmannian. Denote by $T_xM_n$ the tangent space of $M_n$ at $x$ and define
$$\widetilde M_n=\big\{ (x, H) \in M_n \times G(k-p,k) \mbox{, } H \mbox{
  parallel to } T_x M_n \big\}.$$
The volume of $\widetilde M_n$ is called {\it the curvature} of $M_n$, see \cite{Di}. 
We have the following result, see \cite{Det6}.

\begin{theorem} \label{criterium}
Let $M_n$ be a sequence of submanifolds of pure dimension $k-p$ in an open set $X$ of $\C^k$. 
Let $v_n$ be the volume of $M_n$ and $\widetilde v_n$ its curvature. 
Assume that $v_n$ and $\widetilde v_n$ are finite and that the sequence 
$T_n=v_n^{-1}[M_n]$ converge to a current $T$ (this property is always true for some subsequences). 
If $\widetilde v_n=O(v_n)$, then $T$ is woven.
\end{theorem}

Here is another criterium that we will use in the proof of the first main result.

\begin{proposition} \label{criterium_bis}
Let $X$ and $Y$ be two complex manifolds of dimensions $k,l$ respectively. Let $T$ be a positive $(p,p)$-current on $X$ and $S$ is a positive $(q,q)$-current on $Y$. Assume that $T\otimes S$ is woven in $X\times Y$. Then $T$ and $S$ are woven in $X$ and $Y$ respectively.
\end{proposition}
\proof
We will show that $T$ is woven. The same proof works for $S$. Since the problem is local, we can assume that $X$ and $Y$ are bounded open domains in $\C^k$ and $\C^l$ respectively. 
Denote by $\Pi_X$ and $\Pi_Y$ the projections from $X\times Y$ to $X$ and $Y$. Choose a linear projection $\pi:\C^l\to \C^{l-q}$ and a positive smooth form $\Omega$ of maximal degree with compact support in $Y$ 
such that $m:=S\wedge \pi^*(\Omega)$ is a non-zero positive measure. Multiplying $\Omega$ with a constant allows us to assume that $m$ is a probability measure.  We deduce that $(\Pi_X)_*(T\otimes m)=T$.

Define $R:=T\otimes S$ and $\Pi:=\pi\circ\Pi_Y$. Observe that 
$$T\otimes m =(T\otimes S)\wedge \Pi^*(\Omega)=R\wedge \Pi^*(\Omega).$$
We will show that $R\wedge \Pi^*(\Omega)$ is woven. 
By definition of woven current, this property together with the identity $(\Pi_X)_*(T\otimes m)=T$ will imply that $T$ is also woven. 

Since $R$ is woven, it is enough to consider the case where $R$ is the current of integration on a connected manifold $\Gamma$ of dimension $k-p+l-q$. If $\dim \Pi(\Gamma)< l-q$, then $[\Gamma]\wedge \Pi^*(\Omega)=0$. So, we can assume that $\dim\Pi(\Gamma)=l-q$. 
Let $\Sigma$ denote the set of critical values of $\Pi_{|\Gamma}$. Then, for $a\not\in\Sigma$, the intersection 
$\Gamma\cap\Pi^{-1}(a)$ is transversal and is either empty or a smooth manifold of dimension $k-p$.
By Bertini's theorem,  $\Pi^{-1}(\Sigma)\cap\Gamma$ has zero volume in $\Gamma$. The Fubini's theorem (see \cite{Ch} p.334) implies that:
$$[\Gamma]\wedge \Pi^*(\Omega)=\int [\Gamma\cap \Pi^{-1}(a)] d\nu(a),$$
where $\nu$ is the positive measure defined by $\Omega$. Clearly, $[\Gamma]\wedge \Pi^*(\Omega)$ is woven. This completes the proof.
\endproof

We now give a criterion in order to prove that a woven $(1,1)$-current is laminar. The following proposition was independently obtained by Dujardin \cite{Duj}.

\begin{proposition} \label{prop_laminar}
Let $T$ be a woven positive closed $(1,1)$-current of a complex manifold $X$. Assume that 
the local potentials of $T$ are integrable with respect to $T$ (so, the wedge-product $T \wedge T$ is well-defined). If $T \wedge T=0$, then $T$ is laminar.
\end{proposition}

\proof
The problem is local. So, we can assume that $X$ is a domain in $\C^k$. 
By Proposition \ref{woven_def}, we can write $T= \sum_{i=0}^{+ \infty} T_i$, where $T_i$ are nice woven currents 
on $\Qc_{2^{-i}}$ which vanish outside $X$. Denote by $\nu_i$ the measure associated with $T_i$ which is defined on the space of nice hypersurfaces in the components of $\Qc_{2^{-i}}$.

Assume that $T$ is not laminar. Then we can find $\Gamma_0$ in the support of some $\nu_i$ and $\Gamma'_0$ in the support of some $\nu_j$ such that $\Gamma_0 \cap \Gamma'_0$ is a non-empty subvariety of dimension $k-2$. Denote by $D$ (resp. $D'$) the $2^{-i}$-cube (resp. $2^{-j}$-cube) which contains $\Gamma_0$ (resp. $\Gamma'_0$).
Denote also by $V_i^{\epsilon}$ (resp. $V_j^\epsilon$) the set of nice hypersurfaces of $D$ (resp. $D'$) whose distance to $\Gamma_0$ (resp. $\Gamma_0'$) is less than $\epsilon$.
If $\epsilon$ is small enough then for $\Gamma \in V_i^{\epsilon}$ and  $\Gamma' \in V_j^{\epsilon}$ the intersection $\Gamma\cap\Gamma'$ is a non-empty variety of dimension $k-2$ (by Hurwitz's theorem).

Since $\Gamma_0$ is in the support of $\nu_i$ and $\Gamma'_0$ is in the support of $\nu_j$, we have
$\nu_i(V_i^{\epsilon})>0$ and $\nu_j(V_j^{\epsilon})>0$.
Define
 $S:= \int_{V_i^{\epsilon}} [\Gamma] d \nu_i(\Gamma)$ and $S'=\int_{V_j^{\epsilon}} [\Gamma'] d \nu_j(\Gamma')$. Then the geometric intersection
$$S\wedge_g S':=\int_{V_i^{\epsilon}} \int_{V_j^{\epsilon}} [\Gamma \cap \Gamma'] d \nu_i(\Gamma) d\nu_j(\Gamma')$$
is a positive closed $(2,2)$-current on $D\cap D'$ with positive mass. Here, the intersection $\Gamma \cap \Gamma'$ is counted with multiplicity, see e.g. \cite[Ch. III, Prop. 4.12]{Demailly}. 

Now $S\wedge S'$ is well defined in $D\cap D'$ because $S\leq T$ and $S'\leq T$ and $T \wedge T$ is well defined. It is not difficult to show that $S\wedge_g S'= S\wedge S'$, see
\cite[Prop. 2.6]{DDG}. This is a contradiction because $S \wedge S'\leq T \wedge T=0$.
\endproof

\section{Woven structure of Green currents} \label{green_current}

Consider now a holomorphic automorphism on a compact K\"ahler manifold $(X,\omega)$ of dimension $k$. 
Fix a norm on the Hodge cohomology ring $\oplus H^*(X,\C)$. We will use
the following result in order to obtain the first main theorem.

\begin{proposition} \label{criterium_global}
Let $\lambda_n$ be the norm of the operator $(f^n)^*$ on $\oplus H^*(X,\C)$. Let $M$ be a submanifold of codimension $p$ of $X$. Then all limit values of the sequence  
$\lambda_n^{-1} (f^n)^*[M]$ are woven currents.
\end{proposition}
\proof
Since $X$ is a compact K\"ahler manifold, the mass  $\|T\|:=\langle T,\omega^{k-p}\rangle$ of a positive closed $(p,p)$-current depends only on its cohomology class in $\oplus H^*(X,\C)$.
Recall also that when $T$ is given by integration on a manifold $M$, by Wirtinger's theorem, the mass of $T$ is 
$(k-p)!$ times the volume of $M$. By definition of $\lambda_n$, the mass of 
$\lambda_n^{-1} (f^n)^*[M]$  is bounded uniformly on $n$. In particular, the limit values of $\lambda_n^{-1} (f^n)^*[M]$ are positive closed $(p,p)$-currents.

Consider a limit $T$ of a sequence $\lambda_{n_i}^{-1} (f^{n_i})^*[M]$. We have to show that $T$ is woven. For this purpose, we will apply Theorem \ref{criterium} for the restriction of $T$ to charts of $X$. 
It is enough to consider the case where $T\not=0$. This property is equivalent to the fact that
the volume of 
$M_{n_i}:=f^{-n_i}(M)$ increases like  $\lambda_{n_i}$ when $i\to\infty$.

Denote by $TX$ the complex tangent bundle of $X$ and  
define  $E:=\bigwedge^{k-p} TX$ the bundle of holomorphic tangent $(k-p)$-vectors.
Consider also the projectivization $\P(E)$ of $E$. The canonical projection 
$\pi:\P(E)\to X$ defines a holomorphic fibration whose fibers are isomorphic to the projective space of dimension $r-1$ where $r:=\rank(E)$. 

If $a$ is a point in $M$, the tangent space of $M$ at $a$ is defined by a non-zero holomorphic tangent $(k-p)$-vector $v$ which is unique up to a multiplicative constant. So, we can associate to $a$ a point $\widetilde a=(a,[v])$ in  $\pi^{-1}(a)$. When $a$ varies in $M$, the point $\widetilde a$ describes a submanifold $\widetilde M$ of dimension $k-p$ in $\P(E)$.  
On a neighbourhood $U$ of $a$, we can identify $TX$ with the trivial vector bundle  
$U\times \C^k$ and $\P(E)$ to the product $U\times \P^{r-1}$. The Grassmannian  $G(k-p,k)$ is canonically identified to a submanifold of  $\P^{r-1}$ and  $\widetilde M$ coincides with the submanifold of  $U\times G(k-p,k)$ introduced in the definition of the curvature of $M$ in Section \ref{lamination}.

The automorphism $f$ lifts canonically to an automorphism on $TX$. Hence, it also lifts canonically  to an automorphism $\widetilde f$ of  $\P(E)$. This map $\widetilde f$ preserves the fibration $\pi:\P(E)\to X$, that is, 
$\pi\circ \widetilde f = f\circ \pi$. 
If $\widetilde M_n$ is associated to $M_n$, we have $\widetilde M_n=\widetilde f^{-n}(\widetilde M)$. 
In order to apply Theorem \ref{criterium}, it is enough to show that
$\volume(\widetilde M_n)=O(\lambda_n)$. 
According to Proposition 3.18 in Voisin \cite{Voisin}, $\P(E)$ is a compact K\"ahler manifold. So, it is enough to verify that  $\|(\widetilde f^n)^*\|=O(\lambda_n)$
on $\oplus H^*(\P(E),\C)$ for a fixed norm on the last space.

Consider the canonical line bundle  $O(1)$ of $\P(E)$. It is invariant under the action induced by $\widetilde f$. If $h$ denotes the Chern class of $O(1)$ in $H^2(\P(E),\C)$, we have $\widetilde f^*(h)=h$. 
Moreover, according to Lemma  7.32 in \cite{Voisin}, the ring $\oplus H^*(\P(E),\C)$ is generated by $h$ and by the sub-ring $\pi^*(\oplus H^*(X,\C))$. In other words, any class in $\oplus H^*(\P(E),\C)$ can be written as a linear combination of classes of the form  $h^m\smile \pi^*(c)$, where $c$ is a class in $\oplus H^*(X,\C)$.
Finally, we deduce from the above discussion that 
$$(\widetilde f^n)^*(h^m\smile \pi^*(c))= h^m\smile\pi^*(f^n)^*(c).$$
The norm of this class increases at most like $\lambda_n$ when $n\to\infty$. It follows that
 $\|(\widetilde f^n)^*\|=O(\lambda_n)$ which completes the proof.
\endproof

In the rest of this section, we give the proof of Theorem \ref{th1}. We assume that the action of $f^*$ on Hodge cohomology is simple and we will use the same notation given in Introduction. Consider the automorphism $F$ on $X \times X$ defined by $F(x,y)=(f(x),f^{-1}(y))$.
By K\"unneth formula \cite{Voisin}, we have 
$$H^{l,l}(X \times X, \C) \simeq  \bigoplus_{p+p'=l\atop
q+q'=l} H^{p,q}(X , \C) \otimes H^{p',q'}(X , \C).$$
Moreover, $F^*=(f^*, f_*)$ preserves this decomposition. 
It was shown in \cite{Di} that the spectral radius of $f^*$ on $H^{p,q}(X , \C) $ is bounded by 
$\sqrt{d_p d_q}$. It follows that the action of $F^*$ on Hodge cohomology is simple.
The dynamical degree  $d_k(F)$ of order $k$ of $F$ is the maximal one and is equal to $d_s^2$. 

The following result together with Propositions \ref{criterium_bis} and \ref{prop_laminar} imply Theorem \ref{th1}. Note that when $T^+$ is of bidegree $(1,1)$, we have $f^*(T^+ \wedge T^+)=d_1^2 T^+ \wedge T^+$ and we
deduce that $T^+ \wedge T^+=0$ since in this case $d_1^2>d_1>d_2$. The same property holds for
$T^-$.

\begin{proposition} \label{diagonal}
Let  $\Delta$ denote the diagonal of  $X \times X$. Then 
$$d_s^{-2n}(F^n)^* [ \Delta ] \to c T^+ \otimes T^{-},$$
where $c > 0$ is a constant. Moreover, $T^+\otimes T^-$ is woven.
\end{proposition}
\proof
We only have to prove the first assertion. The second one is then a consequence of Proposition \ref{criterium_global} applied to $F$ instead of $f$.

Since the action of  $F^*$ on $\oplus H^*(X\times X,\C)$ is simple and $d_s^2$ is its spectral radius, 
the sequence $d_s^{-2n}(F^n)^*$ converges to a rank 1 linear operator on $\oplus H^*(X\times X,\C)$. The image of the limit operator is a complex line. It is generated by the class  $\{T^+\} \otimes \{T^-\}$ since we have $F^*(T^+\otimes T^-)=d_s^2T^+\otimes T^-$. In particular, if $\{\Delta\}$ denotes the class of the current $[ \Delta]$ in $H^{k,k}(X \times X , \C)$, 
$d_s^{-2n}(F^n)^* \{\Delta\}$  converge in $H^{k,k}(X \times X, \C)$ to 
$c\{T^+\} \otimes \{T^-\}$ for some constant $c$. We have $c\geq 0$ since $[\Delta]$ and $T^\pm$ are positive closed currents. We first show that $c\not =0$.

Denote by  $\Pi_1$ and $\Pi_2$ the projections from $X \times X$ on its factors. On one hand, the integral 
$$\big\langle d_s^{-2n}(F^n)^* [ \Delta ], \Pi_1^* \omega^{k-s} \wedge \Pi_2^* \omega^s \big\rangle$$
(which can be computed cohomologically)  converges to 
$$c\langle T^+ \otimes T^- , \Pi_1^* \omega^{k-s} \wedge \Pi_2^* \omega^s \rangle.$$ 
On the other hand, the same integral is equal to 
\begin{equation*}
\begin{split}
d_s^{-2n} \big\langle [ \Delta ] , (F^n)_*( \Pi_1^* \omega^{k-s} \wedge \Pi_2^* \omega^s )\big\rangle&=d_s^{-2n} \big\langle [ \Delta ] , \Pi_1^* (f^n)_* \omega^{k-s} \wedge \Pi_2^* (f^n)^*\omega^s \big\rangle\\
&=d_s^{-2n} \big\langle (f^n)_* \omega^{k-s} , (f^n)^*\omega^s \big\rangle.
\end{split}
\end{equation*}
The last equality is obtained using that $\Pi_1$ and $\Pi_2$ are equal on $\Delta$.
Finally, the last expression is equal to
$$d_s^{-2n} \big\langle \omega^{k-s} , (f^{2n})^*\omega^s \big\rangle$$
which converges to $ \langle \omega^{k-s} , T^+ \rangle =\|T^+\|\neq 0$. It follows that 
$c \neq 0$.

Now, according to Corollary  4.3.4 in \cite{DS1}, $d_s^{-2n}(F^n)^* [ \Delta ]$ converges to a Green $(k,k)$-current in the cohomology class  $c\{T^+\}\otimes \{T^-\}$. 
It is follows from Theorem 4.3.1 in  \cite{DS1}, that this Green current is the unique positive closed current in its cohomology class. Therefore, it is equal to $cT^+\otimes T^-$. We conclude that  $d_s^{-2n}(F^n)^* [\Delta]$ converges to  $cT^+\otimes T^-$. 
\endproof

\section{Equidistribution towards Green currents} \label{equidistribution}

In this section, we will give some results on the convergence towards the main Green currents of $f$ which will allow us to prove our second main result.

Let $f:X\to X$ be a holomorphic automorphism with simple action on Hodge cohomology as above. Let $\lambda_{p,n}$ denote the norm of $(f^n)^*$ on $H^{p,p}(X,\C)$. Choose a basis on $H^{p,p}(X,\C)$ such that $f^*$ has the Jordan form. We see that $\lambda_{p,n}$ is equivalent to $n^{m_p}d_p^n$ for some positive integer $m_p$ depending on $p$ ($m_p+1$ is the size of a Jordan block). Recall that $d_s$ is the maximal dynamical degree and since the action of $f^*$ on the Hodge cohomology is simple, we have $m_s=1$. 
We also have $\lambda_{0,n}=\lambda_{k,n}=d_0=d_k=1$, 
$d_{p-1}<d_p$ for $p\leq s$ and $d_p>d_{p+1}$ for $p\geq s$. In particular, we have $\lambda_{p-1,n}\lesssim \lambda_{p,n}$ for $p\leq s$. 

\begin{proposition} \label{mass_estimate}
Let $S$ be a positive closed current of bidegree $(k-p,k-p)$ on an open subset $U$ of $X$ with $p \leq s$. 
Let $\chi$ be a positive smooth function with compact support on $U$. Then there is a constant $c>0$ such that for all $r\leq p$ and $n,m\in\N$ we have
$$\|(f^n)_*(\chi S)\wedge (f^m)^*(\omega^r)\|\leq c\lambda_{p,n}\lambda_{r,m}.$$
\end{proposition}

\proof
We first prove the proposition for $r=0$ by induction on $p$. 
The property holds for $p=0$ because in this case  $\chi S$ and $S_n=(f^n)_*(\chi S)$ are positive measures of the same mass.
Assume now the property for $p-1$. We show it for $p$. 
Denote for simplicity $S':=\chi S$. If $\Omega$ is a smooth closed form of bidegree $(p,p)$, we show that 
$$|\langle (f^n)_*(S'),\Omega\rangle|\lesssim \lambda_{p,n}.$$
We then obtain the result by taking $\Omega:=\omega^p$. 

Fix smooth closed $(p,p)$-forms $\alpha_1,\ldots,\alpha_h$ of bidegree $(p,p)$ such that the classes 
$\{\alpha_i\}$ form a basis of $H^{p,p}(X,\C)$ on which $f^*$ has the Jordan form.
Observe that $i\ddbar\chi$ is a closed 
smooth real $(1,1)$-form. So, it can be written as the difference of two positive closed $(1,1)$-forms on $X$. Therefore, we can write $i\ddbar S'=S_1-S_2$ where $S_1, S_2$ are 
two positive closed currents of bidegree $(k-p+1,k-p+1)$ 
on $U$. If $\chi'$ is a smooth positive function with compact support on $U$ and equal to 1 on $\supp(\chi)$, then $\chi'=1$ on $\supp(i\ddbar S')$. So, we can also write $i\ddbar S'=\chi'S_1-\chi'S_2$ and apply  
the induction hypothesis. We have 
$$\|i\ddbar (f^n)_*(S')\|=\|(f^n)_*(i\ddbar S')\|\lesssim \lambda_{p-1,n}\lesssim \lambda_{p,n}.$$

When $\Omega$ is exact, by $\ddbar$-lemma \cite{Voisin}, we can write $\Omega=i\ddbar \Theta$ with $\Theta$ smooth. Since $\Theta$ can be bounded by positive closed forms, it follows that
$$|\langle (f^n)_*(S'),\Omega\rangle|=|\langle i\ddbar (f^n)_*(S'),\Theta\rangle|\lesssim  \lambda_{p,n}.$$
So, the desired estimate holds when $\Omega$ is exact. Subtracting from $\Omega$ an exact form allows us to assume that $\Omega$ is a linear combination of $\alpha_i$. Therefore, it is enough to consider the case where 
$\Omega=\alpha_i$. 

A priori there are several Jordan blocks but we can work with each of them separately.
So, without loss of generality, assume that $\{\alpha_1,\ldots,\alpha_l\}$ corresponds to a Jordan block with eigenvalue $\theta$. 
In particular, we have $f^*\{\alpha_i\}=\theta\{\alpha_i\}+\{\alpha_{i-1}\}$.
We can choose $\alpha_i$ so that 
$f^*(\alpha_i)=\theta\alpha_i+\alpha_{i-1}$ for $2\leq i\leq l$. Indeed, we can first fix $\alpha_l$ and then define the other $\alpha_i$ by induction.
By definition of $d_p$ and $m_p$, we also have $|\theta|\leq d_p$ and when  $|\theta|=d_p$ we should have $l\leq m_p+1$. It suffices to consider the case where $\Omega=\alpha_i$ for $1\leq i\leq l$. 
We show by induction on 
$1\leq i\leq l$ that 
$$|\langle (f^n)_*(S'),\alpha_i\rangle|\lesssim n^{i-1}  d_p^n.$$

Define
$$I_n:=\langle (f^n)_*(S'),\alpha_i\rangle = \langle S', (f^n)^*(\alpha_i)\rangle.$$
Consider the case $i=1$. Since the form
$\beta:=f^*(\alpha_1)-\theta\alpha_1$ is exact, we obtain as above that
$$|I_{j+1}-\theta I_j|= |\langle S',(f^j)^*(\beta)\rangle|= |\langle (f^j)_*(S'),\beta\rangle| \lesssim 
\lambda_{p-1,j}.$$
On the other hand, we have
$$I_n=\sum_{j=0}^{n-1} \theta^{n-1-j} (I_{j+1}-\theta I_j) + \theta^n I_0. $$
It is easy to deduce that $|I_n|\lesssim d_p^n$.
So, the desired inequality is true for $\alpha_1$. 

Assume that the inequality is true for $\alpha_{i-1}$ for 
some $2\leq i\leq l$. We show it for $\alpha_i$. Using that $f^*(\alpha_i)=\theta \alpha_i+\alpha_{i-1}$, we obtain that
$$ \langle S', (f^n)^*(\alpha_i)\rangle = \sum_{j=0}^{n-1} \theta^{n-1-j} \langle S', (f^j)^*(\alpha_{i-1})\rangle +  \theta^n \langle S', \alpha_i\rangle.$$
Therefore,
$$ |\langle S', (f^n)^*(\alpha_i)\rangle|\lesssim  \sum_{j=0}^{n-1} |\theta|^{n-1-j} j^{i-2}d_p^j+|\theta|^n\lesssim 
n^{i-1}d_p^n.$$
This completes the proof of the proposition for $r=0$. 

The general case is also proved by induction on $p$. When $p=0$, we have $r=0$ and the estimate is clearly true. Assume the proposition for $S$ of bidegree $(k-p+1,k-p+1)$. We prove that it is also true for bidegree $(k-p,k-p)$.  If $r=p$, we have $$\|(f^n)_*(\chi S)\wedge (f^m)^*(\omega^p)\|=\langle (f^n)_*(\chi S), (f^m)^*(\omega^p)\rangle 
= \langle (f^{n+m})_*(\chi S), \omega^p\rangle.$$
Using the case $r=0$ above, we can bounded the last integral by a constant times
$\lambda_{p,n+m}$ which is  $\lesssim \lambda_{p,n}\lambda_{p,m}$. So, the proposition is true for $r=p$. 

Assume now that $r\leq p-1$.  Define for simplicity $S'':=(f^n)_*(\chi S)$. 
We show for any smooth closed $(r,r)$-form $\Omega$ that 
$$|\langle S'', (f^m)^*(\Omega)\wedge \omega^{p-r}\rangle|\lesssim \lambda_{p,n}\lambda_{r,m}.$$
We obtain the result by taking $\Omega=\omega^r$. The proof uses the same idea as above. First, we can reduce the problem to the case where $\Omega$ is a form corresponding to 
a Jordan basis of $H^{r,r}(X,\C)$. Then, we follow closely the arguments given in the case $r=0$. The details are left to the reader. Note that our proof is valid for non-invertible maps and we only need to assume that $d_{p-1}<d_p$. 
\endproof

\begin{corollary} \label{limit_closed}
Let $S$ and $\chi$ be as in Proposition \ref{mass_estimate}. If $T$ is a limit value of the sequence $S_n:=\lambda_{p,n}^{-1}(f^n)_*(\chi S)$, then $T$  is a positive closed current.   If $p=s$ then 
$S_n:=d_s^{-n}(f^n)_*(\chi S)$ converge to $c T^-$ where $c \geq 0$ is a constant.
\end{corollary}
\proof
We prove the first assertion. The case $p=0$ is clear since 
$S_n$ is a positive measure and its mass is independent of $n$. So, assume that $p\geq 1$. 
We have seen that $i\ddbar S'$ is the difference of two positive closed currents of bidegree $(k-p+1,k-p+1)$ 
on $U$. 
Applying Proposition \ref{mass_estimate} to $p-1$ instead of $p$ and to $r=0$, we obtain
$$\|i\ddbar S_n\| = \lambda_{p,n}^{-1}\|(f^n)_*(i\ddbar S')\|  \lesssim \lambda_{p,n}^{-1}\lambda_{p-1,n}.$$
Since $d_{p-1} < d_p$, the last expression tends to 0. It follows that $T$ is $\ddbar$-closed. 

Applying the same idea to the map $(f,f)$ on $X\times X$ and to the current $S \otimes S$, 
we obtain that $T \otimes T$ is $\ddbar$-closed (one can easily check that the
degree of order $2s$ of $(f,f)$ is strictly larger than the other dynamical degrees).
Now, we have 
$$0=\partial \overline{\partial} ( T \otimes T) = - \overline{\partial} T \otimes \partial T + \partial T \otimes \overline{\partial} T.$$
By considering the degrees corresponding to each factor of $X\times X$, we see that both terms in the last sum vanish. Hence, $\partial T=0$, $\dbar T=0$ and $T$ is closed.

We now prove the second assertion for $p=s$ and $S_n:=d_s^{-n}(f^n)_*(\chi S)$. 
By Proposition \ref{mass_estimate}, the family of limit values 
of the sequence $S_n$ is a compact set $\Kc$.
Let $S_{n_i}$ be a subsequence which converges to a current $T$. 
We can extract a subsequence of $S_{n_i -1}$ which converges to a current $T_1$. Then, we have 
$d_s^{-1} f_* (T_1)=T$. By induction, we construct a sequence $T_i$ in $\Kc$ with $d_s^{-i} f^{i}_* (T_i) = T$. 

Recall that the action of $f^*$ on the Hodge cohomology is simple. Since the $T_i$'s belong to a compact set and $d_s^{-i} f^{i}_* (T_i) = T$, any limit value of  $d_s^{-i} f^{i}_* \{T_i\}$ belongs to the line generated by
the class $\{T^-\}$. It follows that $\{T\}=c \{T^- \}$ for some $c\in \C$. We have $c\geq 0$ because $T$ and $T^-$ are positive closed currents. 
By uniqueness of Green currents \cite[Th.4.3.1]{DS2}), $T^-$ is the unique positive closed current in its cohomology class. Therefore, we have $T=c  T^-$. It remains to show that the mass $c$ does not depend on the choice of $T$. 

Let $\alpha$ be a smooth $(s,s)$-form in the class $\{T^+\}$. We first show that $\langle S_n , \alpha \rangle$ converge. 
For this purpose, it is enough to verify that $\langle S_n , \alpha \rangle -  \langle S_{n+1} , \alpha \rangle$ decreases exponentially fast. The last expression is equal to
$$\langle S_n , \alpha - d_s^{-1}  f^* (\alpha) \rangle.$$
Since $f^*(T^+)=d_s T^+$, the form $ \alpha - d_s^{-1}  f^* (\alpha)$ is exact. So, we can write it as $i\ddbar \beta$ with $\beta$ a smooth form. Therefore,
$$\langle S_n , \alpha \rangle -  \langle S_{n+1} , \alpha \rangle = \langle i\ddbar S_n , \beta \rangle.$$
Proposition \ref{mass_estimate} implies that the last integral decreases exponentially fast. So,  $\langle S_n , \alpha \rangle$ converge to some constant $c_0$. 

We deduce from this convergence that $\langle T,\alpha\rangle = c_0$. It follows that 
$c\langle T^-,\alpha\rangle=c_0$. Since $\alpha$ is cohomologous to $T^+$ and $T^+\wedge T^-$ is a probability measure, we obtain that $c=c_0$. So, $c$ is independent of the choice of $T$ and this completes the proof.
\endproof

\begin{remark} \label{limit} \rm
Assume that $m_p=0$ and that $S_n\to 0$.
Since, $\lambda_{p,n}^{-1}\lambda_{p-1,n}$ decreases to 0 exponentially fast,  
we deduce from the estimates in the above proofs that the mass of $S_n$ on a compact 
set of $U$ decreases to 0 exponentially fast.
\end{remark}

Corollary \ref{limit_closed} allows us to define the intersection $S\wedge T^+$ of $T^+$ with a positive closed $(k-s,k-s)$-current $S$ on $U$. The intersection is a positive measure given by
$$\langle S \wedge T^+,\chi\rangle:=\lim_{n \rightarrow + \infty} \langle\chi S, d_s^{-n}(f^n)^*(\omega^s)\rangle=\lim_{n \rightarrow + \infty} \langle d_s^{-n}(f^n)_*(\chi S),\omega^s \rangle$$
for $\chi$ smooth with compact support in $U$. 
We have the following result.

\begin{proposition} \label{conv}
Let $S$ be a current of bidegree $(k-s,k-s)$ and $\chi$ a function as above. Then $d_s^{-n} (f^n)_*(\chi S)\wedge T^+$ converge to $c \mu = c T^+ \wedge T^-$, where $c \geq 0$ is the constant such that $d_s^{-n} (f^n)_*(\chi S)\to cT^-$. 
\end{proposition}
\proof
Fix a smooth function $\varphi$ on $X$ and consider
 $$S_n:=d_s^{-n} (f^n)_*(\chi S)= d_s^{-n}(\chi\circ f^{-n})(f^n)_*(S)
 \quad \mbox{and}\quad a_n:=\langle  S_n\wedge T^+,\varphi\rangle.$$
We have $S_n\to cT^-$ and 
$$a_n= \lim_{N \rightarrow \infty} \big\langle \varphi  S_n , d_s^{-N}(f^N)^*(\omega^s) \big\rangle=\lim_{N \rightarrow \infty} \big\langle \varphi  S_n , L^N(\omega^s) \big\rangle,$$
where we denote for simplicity $L:=d_s^{-1}f^*$. 
We want to show that $a_n\to c\langle \mu,\varphi\rangle$. 

We have $T^+=\lim L^N(\omega^s)$.
Fix a smooth closed $(s,s)$-form $\alpha$ in the class $\{T^+\}$. We can complete $\{\alpha\}$ in order to obtain a Jordan basis for $H^{s,s}(X,\C)$. The class $\{\alpha\}$ corresponds to the eigenvalue of maximal modulus $d_s$ of $f^*$ which is a simple eigenvalue. Arguing as in Proposition \ref{mass_estimate}, we obtain that 
$L^N(\omega^s)-L^N(\alpha)\to 0$ and hence $L^N(\alpha)\to T^+$. So,
we can write 
$$T^+=\alpha+\sum_{N=0}^{+ \infty} L^N(L\alpha-\alpha).$$
We also obtain in the same way that 
$$a_n= \lim_{N \rightarrow \infty} \big\langle \varphi  S_n , L^N(\alpha) \big\rangle.$$
This allows us to write
\begin{eqnarray*}
a_n &=& \lim_{N \rightarrow \infty} \Big\langle \varphi  S_n, \alpha+\sum_{m=0}^{N-1} L^m(L\alpha-\alpha)  \Big\rangle\\
&= & \langle \varphi S_n,\alpha\rangle +\sum_{m=0}^{+ \infty} \big\langle \varphi S_n,L^m(L\alpha-\alpha)\big\rangle.
\end{eqnarray*}

We claim that
$$|\langle \varphi S_n,L^m(L\alpha-\alpha)\rangle|\lesssim\gamma^m$$
for some constant $\gamma<1$ independent of $n$ and $m$. 
Assume the claim. We first complete the proof. The estimate allows us to take the limit term by term when $n\to\infty$ in the last identity for $a_n$. It follows that $a_n$ converges to
$$\big\langle \varphi cT^-,\alpha\big\rangle + \sum_{m=0}^{+ \infty} \big\langle \varphi cT^- ,  L^m(L\alpha-\alpha)\big\rangle
=c\big\langle \varphi T^-, \lim_{m \rightarrow + \infty}  L^{m+1}(\alpha)\big\rangle = c\langle \mu,\varphi\rangle.$$
For the convergence $T^-\wedge L^{m+1}(\alpha)\to T^-\wedge T^+=\mu$ see \cite{DS1}.

It remains to prove the claim. Since $L(T^+)=T^+$, $L(\alpha)$ is cohomologous to $\alpha$. So, we can
write $L\alpha-\alpha=i\ddbar \beta$ with $\beta$ a smooth $(s-1,s-1)$-form. 
Adding to $\beta$ a large constant times $\omega^{s-1}$ allows us to assume that $\beta$ is positive.
We have 
$$|\langle \varphi S_n,L^m(L\alpha-\alpha)\rangle|=|\langle i\ddbar (\varphi S_n),L^m(\beta)\rangle|.$$
We expand the current with $\ddbar$ and use Cauchy-Schwarz's inequality in order to bound the last expresion. It is bounded by the sum of the following 4 integrals
$$|\langle i\ddbar \varphi\wedge S_n, L^m(\beta)\rangle|,\quad |\langle\varphi i\ddbar S_n,L^m(\beta)\rangle|$$
and
$$|\langle i\dbar \varphi\wedge \partial S_n, L^m(\beta)\rangle|,\quad |\langle i\partial \varphi \wedge \dbar S_n,L^m(\beta)\rangle|.$$

According to Proposition \ref{mass_estimate},
the first two integrals are $\lesssim d_s^{-m}\lambda_{s-1,m}$. So, they are $\lesssim \lambda^m$ for some $\lambda<1$. Choose a smooth positive function $\widetilde\chi$ with compact support in $U$ such that 
$i\partial \chi\wedge \dbar \chi\leq \widetilde \chi\omega$. Define $\widetilde S:=\omega\wedge S$ and
$\widetilde S_n:=d_s^{-n}(f^n)_*(\widetilde \chi\widetilde S)$.
Since $S_n=d_s^{-n}(f^n)_*(\chi S)$, using Cauchy-Schwarz's inequality, we can bounded the two other integrals by
$$|\langle i\partial\varphi\wedge\dbar\varphi\wedge S_n,L^m(\beta)\rangle|^{1/2} |\langle \widetilde S_n, L^m(\beta)\rangle|^{1/2}.$$
Proposition \ref{mass_estimate} implies that the last expression is also $\lesssim d_s^{-m}\lambda_{s-1,m}$. The result follows.
\endproof

\section{Uniqueness of measure of maximal entropy} \label{entropy}

In this section, we prove Theorem \ref{th2}. 
Consider an invariant probability measure $\nu$ of maximal entropy $h_\nu=\log d_s$. 
We want to show that $\nu=\mu$. 
Recall that entropy is an affine function on the convex compact set of invariant probability measures. The measure $\nu$ can be decomposed as an average of ergodic measures. So, we only have to consider the case where $\nu$ is ergodic. 

Denote by $\chi_1\geq \chi_2\geq\cdots\geq \chi_k$ the Lyapounov exponents of $\nu$. According to \cite[Cor. 3]{DT1}, we have
$$\chi_1 \geq \dots \geq \chi_s \geq \frac{1}{2} \log \frac{d_s}{d_{s-1}}  > 0$$
and
$$ 0>  \frac{1}{2} \log \frac{d_{s+1}}{d_s}  \geq \chi_{s+1} \geq \dots \geq \chi_k.$$

This property allows us to use Pesin's theory (see \cite{Pe} and \cite{LS}). There is a measurable $f^{-1}$-invariant partition $\xi^u$ whose fibers are open subsets of the local unstable manifolds associated with $\nu$ such that 
$$h_{\nu}= h_{\nu}(f, \xi^u).$$
Here, the last expression denotes the entropy of $\nu$ relatively to the partition $\xi^u$. 
These unstable manifolds are of dimension $s$.
We now follow the approach of Bedford-Lyubich-Smillie, see \cite[Prop. 3.2]{BLS}. 

For $\nu$ almost every point $x$, denote by  $\xi^u(x)$ the atom of the partition $\xi^u$
which contains $x$. By construction, $\xi^u(x)$ is an open set of a unstable manifold. So, we can choose an open set $U_x$ such that $\xi^u(x)$ is a submanifold of $U_x$. We can assume that $\xi^u(x)$ admits a holomorphic extension to a neighbourhood of $U_x$. So, as we have seen above,  $T^+ \wedge [\xi^u(x)]$ is a well-defined positive measure in $U_x$. Since, $\xi^u(x)$ admits an extension, this measure is of finite mass and we can considered it as a measure on $X$. Moreover, the obtained measures verify 

$$(f^{i})_* (T^+ \wedge [\xi^u(x)])= T^+ \wedge \frac{[f^{i} (\xi^u(x))]}{d_s^{i}}.$$

Here, is a crucial point of the proof.

\begin{proposition}\label{instable}
We have
$$\|[\xi^u(x)] \wedge T^+\| > 0$$
for $\nu$ almost every point $x$.
\end{proposition}

We need the following result due to Newhouse \cite{Ne} which is valid in a more general setting, see also  \cite[Lemma 5.2]{Du}. Choose relatively compact open subsets $\eta^u(x)$ of $\xi^u(x)$ which contains $x$. 

\begin{proposition}
For $\nu$-almost every $x$, we have 
$$\liminf_{n \rightarrow \infty}  \frac{1}{n} \log \volume(f^n(\eta^u(x))) \geq h_\nu.$$
\end{proposition}

\noindent
{\bf Proof of Proposition \ref{instable}.} 
Let $\chi$ be a smooth positive function with compact support in $U_x$ and equal to 1 on $\eta^u(x)$ for $x$ generic with respect to $\nu$. We deduce from the last proposition that
$$\liminf_{n \rightarrow \infty}  \frac{1}{n} \log 	\| (f^n)_* (\chi [ \xi^u(x) ]) \| \geq \log d_s .$$

It follows from Corollary \ref{limit_closed} and Remark \ref{limit} that  $S_n= d_s^{-n} (f^n)_* (\chi [ \xi^u(x) ])$ converges to $cT^-$ with $c>0$. Therefore,
$$\int \chi [\xi^u(x)] \wedge T^+ =\lim_{n\to\infty}\langle \chi [\xi^u(x)], d_s^{-n}(f^n)^*(\omega^s)\rangle 
=\lim_{n\to\infty}\langle S_n, \omega^s\rangle> 0.$$
The proposition follows. \hfill $\square$

\bigskip

\noindent
{\bf End of the proof of Theorem \ref{th2}.} We obtain exactly as 
in \cite[Prop. 3.2]{BLS} the following lemma.

\begin{lemma}
For $\nu$ almost every point $x$, the measure 
$$\eta_x=T^+ \wedge [\xi^u(x)]/ \rho(x), \quad \mbox{where } \rho(x):= \|T^+ \wedge [\xi^u(x)]\|$$ 
is equal to the conditional measure
$\nu_x:=\nu(\cdot | \xi^u(x))$.
\end{lemma}

Now, also following the same lines in  \cite[Th. 3.1]{BLS}, we obtain
$$\frac{1}{n} \sum_{i=0}^{n-1} f^i_*(\nu_x) \rightarrow \nu$$
for $\nu$ almost every point $x$.
Since $\nu_x$ is equal to $\eta_x$, it remains to show that 
$$\frac{1}{n} \sum_{i=0}^{n-1} f^{i}_*( \eta_x) \rightarrow \mu.$$

But 
$$(f^{i})_* \eta_x = \frac{T^+ \wedge [f^{i} (\xi^u(x))]}{\rho(x) d_s^{i} }\cdot$$
Now, it is enough to apply Proposition \ref{conv} in order to obtain the result.
\hfill $\square$

\bigskip\noindent
Henry De Th\'elin, Laboratoire Analyse, G\'eom\'etrie et Applications, UMR 7539,
Institut Galil\'ee, Universit\'e Paris 13, 99 Avenue J.-B. Cl\'ement, 93430 Villetaneuse, France.\\
 {\tt dethelin@math.univ-paris13.fr}

\bigskip\noindent
Tien-Cuong Dinh, UPMC Univ Paris 06, UMR 7586, Institut de Math{\'e}matiques de Jussieu,
4 place Jussieu, F-75005 Paris, France.\\
{\tt  dinh@math.jussieu.fr}, {\tt http://www.math.jussieu.fr/$\sim$dinh}

\end{document}